\documentclass{amsart}
\usepackage[leqno]{amsmath}

\usepackage{amsmath}
\usepackage{amsthm}
\usepackage{amssymb}
\usepackage{amscd}
\usepackage{amsxtra}     % Use various AMS packages
\usepackage{epsfig}
\usepackage{verbatim}
\usepackage[all]{xypic}

\theoremstyle{plain} %default

\newtheorem{thm}{Theorem}[section]
\newtheorem{lem}[thm]{Lemma}
\newtheorem{prop}[thm]{Proposition}

\newtheorem{cor}[thm]{Corollary}

\newtheorem{eg}[thm]{Example}

\newtheorem{conj}[thm]{Conjecture}

\theoremstyle{remark}
\newtheorem*{rmk}{Remark}

\numberwithin{equation}{section}

\newcommand{\PP}{\mathbf{P}}

\newcommand{\CC}{\mathbf{C}}

\newcommand{\ZZ}{\mathbb{Z}}
\newcommand{\QQ}{\mathbb{Q}}

\newcommand{\tensor}{\otimes}

 \DeclareMathOperator{\Tor}{Tor}
 \DeclareMathOperator{\Ext}{Ext}
 \DeclareMathOperator{\Hom}{Hom}
 \DeclareMathOperator{\Torv}{\widehat{Tor}}
 \DeclareMathOperator{\Extv}{\widehat{Ext}}

 \DeclareMathOperator{\Spec}{Spec}
 \DeclareMathOperator{\Sing}{Sing}
 \DeclareMathOperator{\Se}{S}
\DeclareMathOperator{\Reg}{R}
 \DeclareMathOperator{\CH}{CH}
 \DeclareMathOperator{\Proj}{Proj}
 \DeclareMathOperator{\im}{im}

 \DeclareMathOperator{\pd}{pd}

 \DeclareMathOperator{\depth}{depth}
 \DeclareMathOperator{\coker}{coker}

 \DeclareMathOperator{\LC}{H}
 
 \DeclareMathOperator{\syz}{syz}

 \DeclareMathOperator{\h}{H}

\newcommand{\Ann}{\textup{Ann}}

\begin{document}

\bibliographystyle{plain}

\title{Some observations on local and projective hypersurfaces}
%\author{Hai Long Dao}
\author{Hailong Dao}
\address{Department of Mathematics \\ University of Utah \\155 South 1400 East, Salt Lake City,
UT 84112-0090, USA} \email{hdao@math.utah.edu} \maketitle
\begin{abstract}
Let $R$ be a hypersurface in an equicharacteristic or unramified
regular local ring. For a pair of modules $(M,N)$ over $R$ we study
applications of  rigidity of $\Tor^R(M,N)$, based on ideas by
Huneke, Wiegand and Jothilingam. We then focus on the hypersurfaces
with isolated singularity and even dimension, and show that modules
over such rings behave very much like those over regular local
rings. Connections and applications to projective hypersurfaces such
as  intersections of subvarieties and cohomological criteria for
splitting of vector bundles are discussed.
\end{abstract}

\section{Introduction}

The purpose of this note is to continue our investigation of the
rigidity and decent intersection of modules over a local
hypersurface $R$ done in \cite{Da1}. There we showed that the
vanishing of $\theta^R(M,N)$, a function introduced by Hochster in
\cite{Ho1}, implies the rigidity of $\Tor$ for $(M,N)$. We will
apply this to give various new results on modules over
hypersurfaces. We also give supporting evidence for the following
Conjecture made in \cite{Da1}:

\begin{conj}\label{conj}
Let $R$ be an hypersurface  with an isolated singularity. Assume
that $\dim R$ is even. Then $\theta^R(M,N)$ always vanishes.
\end{conj}

This conjecture can be viewed as a consequence of a local version of
the following conjecture by Hartshorne on Chow groups of smooth
projective hypersurfaces :

\begin{conj}\label{conjHar}(R.Hartshorne, \cite{Ha1}, page 142) Let $X$ be a smooth
projective hypersurface in $\mathbb{P}^n_{\CC}$. Then $\CH^i(X)=
\mathbb{Z}$  for $i<{\dim X}/2$.
\end{conj}

Note that the original form of Hartshorne's question states that in
codimension less than  ${\dim X}/2$, a cycle is homologically
equivalent to $0$ if and only if it is rationally equivalent to $0$.
It is not hard to see that Hartshorne's original statement is
equivalent to the version stated above. For a K-theoretic discussion
of \ref{conjHar}, see \cite{Pa} (Conjecture 1.5 and Section 6). For
more discussions on how \ref{conjHar} is related to \ref{conj}, see
\cite{Da1}, Section 3. We remark that this connection actually
motivated this project, as well as the proof of \ref{projhyper}
below.

In Section 2 we review the basic notations and preliminary results.
In Section 3 we use the procedure invented by Huneke and Wiegand in
\cite{HW2} to show that $\Hom_R(M,M)$ rarely has good depth. The use
of $\theta^R(M,N)$ allows us to simplify  and strengthen some
results in \cite{HW2} (see \ref{lcrigid} and \ref{3.2}, also
\ref{coreven1} and \ref{coreven2}).

Section 4 is concerned with two consequences of Conjecture
\ref{conj}:

\begin{conj}\label{conj1}
Let $R$ be an admissible hypersurface (meaning $\hat R$ is a
quotient of an unramified or equicharacteristic regular local ring
by a nonzero element) with an isolated singularity . Assume that
$\dim R$ is even. Let $M,N$ be $R$-modules such that
$\Tor_i^R(M,N)=0$ for some $i$. Then $\Tor_j^R(M,N)=0$ for all
$j\geq i$.
\end{conj}

\begin{conj}\label{conj2}
Let $R$ be an admissible hypersurface with an isolated singularity.
Assume that $\dim R$ is even. Let $M,N$ be $R$-modules such that $\
l(M\tensor_RN)<\infty$. Then $\dim M +\dim N \leq \dim R$.
\end{conj}

We will prove \ref{conj} when $M$ is free on the punctured spectrum
of $R$ and $N=M^*$ (see \ref{even}). We also prove \ref{conj2} when
$R$ is a standard graded local hypersurface and $M,N$ are graded
modules (see \ref{graded}). Many applications follow, such as
generalizations of results on $\Hom_R(M,M)$ over regular local rings
by Auslander and Auslander-Goldman (Corollary \ref{coreven3} and
\ref{coreven4}). In that section, the connection with geometry of
projective hypersurfaces (which goes both ways) will be exploited.
For example, \ref{graded} will be proved using $l$-adic cohomology,
and from \ref{even} we obtain the following (see \ref{lines}):

\begin{cor}
Let $k$ be a field and $n$ an even integer $\geq 4$. Let $X \subset
\PP_k^{n}$ be a nonsingular hypersurface. Let $E$ be a vector bundle
on $X$. If $H^1(X, (E\tensor E^*)(l))=0$ for all $l \in \ZZ$, then
$E$ splits.
\end{cor}

In Section 5, we give some further applications of rigidity. The
first is a simple characterization of maximal Cohen-Macaulay modules
(see \ref{MCM}). The second is a connection between  vanishing of
$\Ext_R^n(M,M)$ and $\pd_RM$, generalizing a result by Jothilingam
(see \ref{Jo1}).

\section{Notation and preliminary results}\label{2}

Unless otherwise specified, all rings are Noetherian, commutative
and local, and all modules are finitely generated. A local ring
$(R,m,k)$ is a \textit{hypersurface} if its completion $\hat R$ has
the form $T/(f)$, where $T$ is a regular local ring and  $f$ is in
the maximal ideal of $T$. We say that  $R$ is admissible (as a
hypersurface) if $T$ is a power series ring over a field or over a
discrete valuation ring.

For a ring $R$ and a non-negative integer $i$, we set $X^i(R) :=
\{p\in \Spec(R)| \dim(R_p)\leq i\}$. We denote by $Y(R)$ the set
$X^{\dim(R)-1}$, the punctured spectrum of $R$. We say that $R$
satisfy the condition $(\Reg_i)$ if $R_p$ is regular for any $p\in
X^i(R)$. We denote by $G(R)$ the Grothendieck group of finitely
generated modules over $R$ and by $\bar{G}(R):= G(R)/[R]$, the
reduced Grothendieck group. Also, we let $\Sing(R) := \{p\in
\Spec(R)| R_p \ \text{is not regular} \}$ be the singular locus of
$R$. For an abelian group $G$, we let $G_{\mathbb{Q}} =
G\tensor_{\mathbb{Z}}\mathbb{Q}$.

Let $M^{*} := \Hom(M,R)$ be the dual of $M$. The module $M$ is
called \textit{reflexive} provided the natural map $M\to M^{**}$ is
an isomorphism. The module $M$ is called \textit{maximal
Cohen-Macaulay} (or sometimes abbreviated as MCM) if $\depth_RM =
\dim R$.

For a non-negative integer $n$, $M$ is said to satisfy $(\Se_n)$ if
:
$$ \depth_{R_p}M_p \geq \min\{n,\dim(R_p)\} \ \forall p\in \Spec(R)$$
(The depth of the $0$ module is set to be $\infty$.).

A pair of $R$-modules $(M,N)$ is called {\textit{rigid}} if for any
integer $i\geq0$, $\Tor_i^R(M,N)=0$ implies $\Tor_j^R(M,N)=0$ for
all $j\geq i$. Moreover, $M$ is \textit{rigid} if for all $N$, the
pair $(M,N)$ is rigid.

One defines the finite length index of the pair $(M,N)$ as :
$$ f_R(M,N) := \min\{ i |\ l(\Tor_j^R(M,N))<\infty \ \text{for $j \geq i$} \} $$

\textbf{The function $\theta^R(M,N)$}\\

Let $R = T/(f)$ be an admissible local hypersurface. The function
$\theta^R(M,N)$ was introduced by Hochster ([Ho1]) for any pair of
finitely generated modules $M,N$ such that $f_R(M,N)<\infty$  as:
$$ \theta^R(M,N) = l(\Tor_{2e+2}^R(M,N)) - l(\Tor_{2e+1}^R(M,N)) .$$
where $e$ is any integer $\geq d/2$. It is well known (see
\cite{Ei}) that $\Tor^R(M,N)$ is periodic of period at most 2 after
$d+1$ spots, so this function is well-defined. The theta function
satisfies the following properties (see \cite{Ho1}). First, if
$M\tensor_RN$ has finite length, then:
$$\theta^R(M,N) = \chi^T(M,N).$$
Secondly, $\theta^R(M,N)$ is biadditive on short exact sequence,
assuming it is defined. Specifically, for any short exact sequence:
$$0 \to N_1 \to N_2 \to N_3 \to 0$$
and any module $M$ such that $f_R(M,N_i)<\infty$ for all $i=1,2,3$,
we have $\theta^R(M,N_2) = \theta^R(M,N_1) + \theta^R(M,N_3)$.
Similarly, $\theta(M,N)$ is additive on the first variable.

In \cite{Da1}, we show that when $\theta^R(M,N)$ can be defined and
vanishes, then $(M,N)$ is rigid:

\begin{prop}\label{rg1}
Let $R$ be an admissible hypersurface and $M,N$ be $R$-modules such
that $f_R(M,N)<\infty$ (so that $\theta^R(M,N)$ can be defined).
Assume $\theta^R(M,N)=0$. Then $(M,N)$ is rigid.
\end{prop}

The following corollary will be used frequently in this note:

\begin{cor}
Let $R$ be an admissible hypersurface and $M,N$ be $R$-modules such
that:
\begin{enumerate}
\item $\pd_{R_p}M_p<\infty$ for any $p \in Y(R)$, the punctured spectrum of $R$
(in particular, this is always true if $R$ has only isolated
singularity).
\item $[N]=0$ in $\overline G(R)_{\mathbb{Q}}$.
\end{enumerate}
Then $\theta^R(M,N)$ can be defined and equals $0$. Consequently,
$(M,N)$ is rigid.
\end{cor}

\begin{proof}
The condition on $M$ ensures that $\theta^R(M,N)$ can be defined for
all $R$-module $N$. In other words, $\theta^R(M,-)$ gives a
$\ZZ$-linear map from $G(R)$ to $\ZZ$. The conclusions are now
obvious (note that $\theta^R(M,R)=0$).
\end{proof}
\textbf{The Pushforward}\\

Let $R$ be a Gorenstein ring and $M$ a torsion-free (equivalent to
$(\Se_1)$) $R$-module. Consider a short exact sequence : $$0 \to W
\to R^{\lambda} \to M^* \to 0 $$ Here $\lambda $ is the minimal
number of generators for $M^*$. Dualizing this short exact sequence
and noting that $M$ embeds into $M^{**}$ we get an exact sequence:
$$ 0 \to M \to R^{\lambda} \to M_1 \to 0$$
This exact sequence is called the \textit{pushforward} of $M$. The
following proposition is taken from [HJW].

\begin{prop}(\cite{HJW}, 1.6)\label{pushforward}
Let $R,M,M_1$ as above. Then for any $p\in \Spec(R)$:\\
(1) $M_p$ is free if and only if $(M_1)_p$ is free. \\
(2) If $M_p$ is a maximal Cohen-Macaulay $R_p$-module, then so is $(M_1)_p$.\\
(3) $\depth_{R_p}(M_1)_p \geq \depth_{R_p}M_p -1$.\\
(4) If $M$ satisfies $(\Se_k)$, then $M_1$ satisfies $(\Se_{k-1})$.
\end{prop}

\textbf{The depth formula}\\

A result by Huneke and Wiegand showed that when  all the high $\Tor$
modules vanish, the depths of the modules satisfy a remarkable
equation:

\begin{prop}(\cite{HW1}, 2.5) Let $R$ be a local complete intersection. Let $M,N$ be non-zero
finitely generated modules over $R$ such that $\Tor_i^R(M,N)=0$ for
all $i\geq 1$. Then:
$$ \depth(M) + \depth(N) = \depth(R)+\depth(M\tensor_RN) $$

\end{prop}

\section{On depth of $\Hom_R(M,M)$ }

Throughout this section $R$ is a local hypersurface. We will call an
$R$ module $M$ such that $M$ is locally free over the punctured
spectrum of $R$ a {\it vector bundle} over $Y(R)$ (or just vector
bundle). We observe that for a vector bundle $M$ and any $R$-module
$N$, $\theta^R(M,N)$ is always defined. In this section we will show
that for certain modules over admissible hypersurfaces,
$\Hom_R(M,M)$ and $M\tensor_R M^*$ rarely have good depth. We will
follow the same procedure as in \cite{HW2}, but with two essential
additions : we focus on vector bundles from the beginning and we
will exploit the function $\theta^R(M,N)$ rather heavily. These will
allow us to simplify and strengthen some results in \cite{HW2}.

\begin{prop}\label{3.1}
Let $R$ be an admissible hypersurface. Let $M$ be a vector bundle
over $Y(R)$ such that $\depth(M) \geq 1$. Let $N$ be an $R$-module
such that $\theta^R(M,N)=0$ and $\depth(M\tensor_R N)\geq 1$. Then
$\Tor_i^R(M,N)=0$ for all $i>0$.
\end{prop}

\begin{proof}
The assumptions ensure that $M$ is $\Se_1$. Hence we have the
pushforward of $M$:
$$ 0 \to M \to R^{\lambda} \to M_1 \to 0$$
We can tensor with $N$ to get  :
$$0 \to \Tor_1^R(M_1,N) \to M\tensor_RN \to N^{\lambda} \to M_1\tensor_RN \to 0$$
By \ref{pushforward}, $M_1$ is also a vector bundle. So
$l(\Tor_1^R(M_1,N))<\infty$. Since it embeds into a module of
positive depth, $\Tor_1^R(M_1,N)$ must be $0$. Clearly
$\theta^R(M_1,N)$ is defined and equal to:
$$ \theta^R(R^{\lambda},N)-\theta^R(M,N) = 0$$
So by \ref{rg1}, $\Tor_i^R(M_1,N) = 0$ for all $i>0$. Since $M \cong
\syz_1^R(M_1)$, we are done.

\end{proof}

The next result is an analogue of Theorem 2.4 in \cite{HW2}
\begin{thm}\label{lcrigid}
Let $(R,m)$ be an admissible hypersurface of dimension $d$. Let $r$
be an integer such that $0 \leq r <d$. Let $M$ be a vector bundle
over $Y(R)$ such that $\depth(M) \geq r$.  Let $N$ be an $R$-module
satisfying $(\Se_r)$, and assume  $\theta^R(M,N)=0$ and $\LC_{m}^{r}
(M\tensor_RN)=0$. Then  $\depth(M\tensor_RN) \geq r+1$, and if $r>0$
we have $\Tor_i^R(M,N)=0$ for all $i>0$.
\end{thm}

\begin{proof}
We will use induction on $r$. If $r=0$ the conclusion is trivial.
Now we assume $r>0$. Then $N$ satisfies $(\Se_1)$ and we have the
pushforward:
$$ 0 \to N \to R^{\lambda} \to N_1 \to 0$$
Tensoring with $M$ we get:
$$0 \to \Tor_1^R(N_1,M) \to N\tensor_RM \to M^{\lambda} \to N_1\tensor_RM\to 0$$
which we break into two short exact sequences:
$$0 \to \Tor_1^R(N_1,M) \to N\tensor_RM \to C \to 0$$
and:
$$0 \to C \to M^{\lambda} \to N_1\tensor_RM\to 0$$

Since $M$ is a vector bundle, $T = \Tor_1^R(N_1,M)$ has finite
length. In particular $\LC_m^{r+1}(T) =0$. By applying $\LC_m^0(-)$
to the first exact sequence we get $\LC_m^{r}(C)=0$. Now applying
$\LC_m^0(-)$ to the second exact sequence and using
$\LC_m^{r-1}(M)=0$ (because $\depth(M)\geq r$) and
$\LC_m^r(M\tensor_R N)=0$ we get $\LC_m^{r-1}(N_1\tensor_R M)=0$. We
need to check the other inductive assumptions for $N_1$. Clearly,
$N_1$ is $(\Se_{r-1})$ by \ref{pushforward}. Also, $\theta^R(M,N_1)$
is defined and equal to:
$$ \theta^R(M,R^{\lambda})-\theta^R(M,N) = 0$$
So by induction we have $\depth(M\tensor_R N_1)\geq r$. Then by
\ref{3.1} $\Tor_i^R(N_1,M)=0$ for all $i>0$, so we have the last
assertion and an exact sequence:
$$0 \to N\tensor_RM \to M^{\lambda} \to N_1\tensor_RM\to 0$$
Therefore  $\depth(M\tensor_RN) \geq r$, and since $\LC_{m}^{r}
(M\tensor_RN)=0$ it follows that $\depth(M\tensor_RN) \geq r+1$ .
\end{proof}

\begin{prop}\label{3.2}
Let $R$ be a local hypersurface  and  $M$ be an $R$-module. Assume
that $\depth(M \tensor_RM^{*})\geq 2$ and $\Tor_i^R(M,M^{*})=0$ for
all $i>0$. Then $M$ is free.
\end{prop}

\begin{proof}
By the depth formula we get:
$$ \depth(M)+ \depth(M^{*}) \geq \dim R +2$$
so by (\cite{Va}, 3.3.16) we must have $\depth(M) = \depth(M^{*})
=\dim R$. On the other hand, the vanishing of all
$\Tor_i^R(M,M^{*})$ forces one of the modules to have finite
projective dimension (see (\cite{HW2}, 1.9). Either way, $M$ must be
free (since $M,M^*$ are maximal Cohen-Macaulay).
\end{proof}

\begin{thm}\label{3.3}
Let $R$ be an admissible hypersurface satisfying condition
$(\Reg_2)$. Let $M$ be a reflexive $R$-module such that one of the
following is satisfied:
\\1) $\theta^R(M,M^*)$ is defined and
equals $0$.
\\1') $[M]=0$ in $\overline G(R)_{\mathbb{Q}}$.
\\If $\Hom_R(M,M)$ satisfies $(\Se_3)$ then $M$ is free.
\end{thm}

\begin{proof}
We use induction on $d=\dim R$. If $d\leq 2$ then $R$ is regular,
and $M$, being reflexive, must be free. Assume $d\geq 3$ and (by
induction hypothesis, since all the conditions localize) that $M$ is
free on the punctured spectrum. In other words, $M$ is a vector
bundle. Consider the natural map $\phi : M^{*}\tensor_RM \to
\Hom(M,M)$. Since $M$ is a vector bundle, the kernel and cokernel of
$\phi$ both have finite length. By considering the long exact
sequence of local cohomology from the two short exact sequences:
$$ 0 \to \ker(\phi) \to M^{*}\tensor_RM \to \im(\phi) \to 0$$
$$ 0 \to \im(\phi) \to \Hom(M,M) \to \coker(\phi) \to 0$$
and using $\LC_m^i(\Hom(M,M))=0$ for $i<3$ (because $\Hom(M,M)$ is
$(\Se_3)$ and $d\geq 3$) we can deduce that
$\LC_m^2(M^{*}\tensor_RM)=0$. Now $M^{*}$ is also a  vector bundle,
so $\theta^R(-,M^*) $ is always defined. Therefore if $[M]=0$ in
$\overline G(R)_{\mathbb{Q}}$ then $\theta^R(M,M^{*})=0$. So in both
cases 1) and 1'), $\theta^R(M,M^{*})=0$ and Theorem \ref{lcrigid}
implies $\depth(M\tensor_RM^*)\geq 3$ as well as
$\Tor_i^R(M,M^{*})=0$ for all $i>0$. The result now follows from
Proposition \ref{3.2}.
\end{proof}

\begin{eg}\label{mainex}
This will be our main example throughout this note. Let
$R=k[[x,y,u,v]]/(xu-yv)$ with $m=(x,y,u,v)$ and $M = (x,y)$. We
claim that  $M^* \cong (x,v)$. Any $R$-linear map $\phi$ from $M$ to
$R$ is determined by $\phi(x)$. Hence $M^*$ is isomorphic to $\{a
\in R | \ ya \in xR \}$, which is easily seen to be $(x,v)$. So
$M\tensor_RM^*\cong (x^2,xy,xv,yv) \cong (x^2,xy,xv,xu) \cong m$
(see Example 1.8 of \cite{HW2}). Using the long exact sequence of
local cohomology for the sequence $0 \to m \to R \to k \to 0$ we get
$\LC_m^1(M\tensor_RM^*) = \LC_m^1(m) = k$ and
$\LC_m^2(M\tensor_RM^*) = \LC_m^2(m) = 0 $. Clearly
$\depth(M\tensor_RM^*)= \depth(m) = 1$. Obviously, $M$ is not free.

Note that both $M$ and $M^*$ are maximal Cohen-Macaulay modules over
$R$. Since $R$ has an isolated singularity, they are also vector
bundles over $Y(R)$. It can be easily computed that
$\Tor_1^R(M,M^*)=k$, $\Tor_2^R(M,M^*)=0$ and $\theta^R(M,M^*)=-1$.

Consider Proposition \ref{3.1} and Theorem \ref{lcrigid}. Let
$N=M^*$. Then the example shows that the condition $\theta^R(M,N)=0$
can not be dropped.

Consider Theorem \ref{3.3}. Note that $\overline G(R)_{\mathbb{Q}} =
\QQ[M] \cong \QQ$. The example shows that the condition $M=0$ in
$\overline G(R)_{\mathbb{Q}}$ can not be dropped either.

\end{eg}

\section{Isolated hypersurface singularities of even dimensions}
In this section we will show some supporting evidence for Conjecture
\ref{conj}. Our results indicate that modules over isolated
hypersurface singularities of even dimensions behave very similarly
to those over regular local rings. We first prove:

\begin{thm}\label{even}
Let $R$ be a hypersurface with isolated singularity. Assume that
$\dim R$ is even. Then for any vector bundle $M$ on $Y(R)$,
$\theta^R(M,M^{*})=0$.
\end{thm}

We need to review the concept of stable (co)homology (for more
details, see \cite{Bu} or \cite{AB}). Let $R$ be a Noetherian ring,
a {\it {complete resolution}} of an $R$-module $M$ is a complex
$\bold T$ such that $\h_n(\bold T) = \h^n(\bold T^*)=0$ for all $n
\in \ZZ$ and $\bold T_{\geq r} = \bold P_{\geq r}$ for some
projective resolution $\bold P$ of $M$ and some integer $r$. It is
known that the modules $\h^i(\Hom_R(\bold T,N))$ and $\h_i(\bold
T\tensor_R N)$ are independent of the resolution, and one calls them
$\Extv^i_R(M,N)$ and $\Torv_i^R(M,N)$, respectively.

Before moving on we recall the Local Duality Theorem in our context
and some consequences that will be used. Let $d=\dim R$. Let
$^{\vee}$ denote $\Hom_R(-,E(k))$, the Matlis dual. Then for any
module $M$ we have an isomorphism:
$$\Ext_R^i(M,R) \cong \LC_{m}^{d-i}(M)^{\vee} $$

In particular, if $M$ is maximal Cohen-Macaulay, then
$\Ext_R^i(M,R)=0$ for $i>0$ and if $l(M)<\infty$ then
$\Ext_R^i(M,R)=0$ for $i\neq d$. Also, if $M$ is maximal
Cohen-Macaulay then so is $M^*$, since by dualizing a free
resolution of $M$  one can see that $M^*$ is a syzygy of infinitely
high order.

We will  need the standard and easy results below, reproved for the
reader's convenience (some of these could be found in \cite{AB}, but
we could not find a full reference):

\begin{lem}\label{facts}
Let $R$ be a local hypersurface. Let $M,N$ be $R$-modules such that
$M$ is maximal Cohen-Macaulay (MCM).Then we have the following
isomorphism:\\
(1) $\Tor_i^R(M,N) \cong \Torv_i^R(M,N)$ for all $i>0$.\\
(2) $\Ext^i_R(M,N) \cong \Extv^i_R(M,N)$ for all $i>0$.\\
(3) $\Torv_i^R(M,N) \cong \Torv_{i+2}^R(M,N)$ for all $i$.\\
(4) $\Extv^i_R(M,N) \cong \Extv^{i+2}_R(M,N)$ for all $i$.\\
(5) $\Torv_i^R(M,N) \cong \Extv^{-i-1}_R(M^*,N)$ for all $i$.\\
(6) $\Tor_i^R(M,N) \cong \Extv^{i+1}_R(M^*,N)$ for all $i>0$.
\end{lem}

\begin{proof}
Let $\bold F : ... \to F_1 \to F_0 \to M^* \to 0 $ be a free
resolution of $M^*$. Since $M$ is MCM we have $M^{**}=M$ and
$\Ext^n(M^*,R)=0$ for all $n> 0$. So dualizing $\bold F$ we get an
exact sequence:
$$\bold F^{*} : \ \ 0 \to M \to F_0^{*} \to F_1^{*} \to ...    $$
Now let $\bold G : ... \to G_1 \to G_0 \to M \to 0 $ be a minimal
free resolution of $M$. We can splice $\bold G$ and $\bold F^{*}$
together to get an exact sequence:
$$\bold T : ... \to G_1 \to G_0 \to F_0^{*} \to F_1^{*} \to ...    $$
It is obvious that $\bold T$ is a complete resolution of $M$. That
proves (1) and (2). Since $M,M^*$ are MCM and $R$ is a hypersurface,
$\bold G, \bold F$  are periodic of period at most $2$. In
particular, $\syz_{2n}^RM^*\cong M^*$ for any $n>0$.  Let $M_i = \im
(F_{i-1}^* \to F_i^*)$. Then $M_{2n} \cong (\syz_{2n}^RM^*)^* \cong
M^{**} \cong M$. Thus for any $n>0$:
$$... \to G_1 \to G_0 \to F_0^{*} \to F_1^{*} \to ... \to F_{2n-1}^{*}  \to 0 $$
is a free resolution of $M$. Hence (3) and (4) follow. From the
construction, $\bold T^*[-1]$ is a complete resolution of $M^*$, and
the canonical isomorphism $\Hom_R(\bold T^*,N) \cong(\bold
T\tensor_R N)$ gives (5). Finally, (6) follows from combining (1),
(5) and (4).
\end{proof}

We will also need the following result by Buchweitz:

\begin{prop}(\cite{Bu},10.3.3)\label{Bu} Let $R$ be a local hypersurface with
isolated singularity such that $d= \dim R$ is even. Then for any two
MCM $R$-modules $M,N$, and integers $i,j$ such that $i-j$ is odd :
$$l(\Extv^i_R(M,N)) = l(\Extv^j_R(M^*,N^*))$$
\end{prop}

\begin{proof}
Since \cite{Bu} is not available publicly, and in any case the
assertion was derived there from very general theory, we will
summarize a self-contained proof. Let $T$ be the complete resolution
of $M$ constructed in the proof of \ref{facts} and $I$ be an
(finite) injective resolution of $R$ over $R$. One can see that the
total complexes associated to the double complexes
$\Hom_R(\Hom_R(T,N),I)$ and $\Hom_R(T^*,\Hom(N,I))$ are isomorphic.
Note that $T^*[-1]$ is a complete resolution of $M^*$. We get two
spectral sequences converging to the same limit:
$$^1E_2^{i,j} = \Ext^i_R(\Extv^{-j}_R(M,N),R)\Rightarrow H^{i+j} $$
and
$$^2E_2^{i,j} = \Extv^{i-1}_R(M^*, \Ext^{j}_R(N,R)\Rightarrow H^{i+j} $$
As $N$ is MCM, $\Ext^j_R(N,R)=0$ for $j>0$, so the second sequence
collapses, leaving $H^{i+j} \cong \Extv^{i+j-1}_R(M^*,N^*)$. On the
other hand, $M$ is free on the punctured spectrum, so
$l(\Extv^{-j}_R(M,N))$ is finite. Thus $^1E_2^{i,j} =0 $ unless
$i=d$, which gives $H^{d+j} \cong \Ext^d_R(\Extv^{-j}_R(M,N),R)
\cong \Extv^{-j}_R(M,N)^{\vee}$. Putting everything together we get:
$$\Extv^{-j}_R(M,N)^{\vee} \cong  \Extv^{d+j-1}_R(M^*,N^*) $$
Since $d$ is even and taking Matlis dual preserves length, part (4)
of \ref{facts} finishes our proof.
\end{proof}

Now we can prove \ref{even}:
\begin{proof} (of Theorem \ref{even})
First, we will prove the theorem for the case when $M$ is MCM. In
this situation, for any integer $i>0$, by (6) of \ref{facts}  we
have:
$$\Tor_i^R(M,M^*) \cong \Extv^{i+1}_R(M^*,M^*)$$
and
$$\Tor_{i+1}^R(M^*,M) \cong \Extv^{i+2}_R(M,M)$$
Hence \ref{Bu} (applied for $N=M$) shows that $\theta^R(M,M^*)=0$.

Now, assume $M$ is a vector bundle. Let $K = \syz^R_1(M)$. We want
to prove that $\theta^R(M,M^*) = \theta^R(K,K^*)$. Dualizing the
short exact sequence: $0\to K\to F\to M\to 0$ we get an exact
sequence:
$$ 0 \to M^* \to F^* \to K^* \to \Ext^1_R(M,R) \to 0 $$
So $M^* + K^* = \Ext^1_R(M,R)$ in $\overline G (R)$. Note that
$\Ext^1_R(M,R)$ has finite length because $M$ is free on the
punctured spectrum of $R$. We claim that any module of finite length
is equal to $0$ in $\overline G (R)_{\mathbb{Q}}$. Since any module
of finite length is a multiple of $[k]$, it is enough to prove the
claim for one finite length module. If $\dim R=0$ then there is
nothing to prove. If $\dim R>0$ then pick a prime $p$ such that
$\dim R/p=1$ and a non-unit $x\notin p$. The exact sequence $0\to
R/p \to R/p \to R/(p+x) \to 0$ shows that $[R/(p+x)]=0$, which is
all we need. So $[M^*] =-[K^*]$ and $\theta^R(M^*,-) =
-\theta^R(K^*,-)$. Since $\theta^R(M,-) = -\theta^R(K,-)$, we have
$\theta^R(M,M^*) = \theta^R(K,K^*)$.

Repeating the equality above we get $\theta^R(M,M^*) =
\theta^R(K,K^*)$ when $K = \syz^R_n(M)$ for any $n>0$. But for
$n\gg0$ $K$ is an MCM $R$-module, so $\theta^R(K,K^*)=0$.
\end{proof}

\begin{cor}\label{coreven1}
Let $R$ be an admissible hypersurface such that $R$ has an isolated
singularity and $\dim R$ is an even number greater than $2$. Let $M$
be a reflexive $R$-module. If $\Hom_R(M,M)$ satisfies $(\Se_3)$ then
$M$ is free.
\end{cor}

\begin{proof}
By localizing on $p \in Y(R)$ and using Theorem \ref{3.3} (note that
a regular local ring can also be considered a hypersurface) we may
assume that $M$ is a vector bundle. The result then follows from
\ref{3.3} and \ref{even}.
\end{proof}

\begin{cor}\label{coreven2}
Let $R$ be an admissible hypersurface such that $R$ has isolated
singularity and $\dim R > 2$ and is even. Let $M$ be a vector bundle
over $Y(R)$. Assume that $M$ is reflexive and
$\LC_m^2(M\tensor_RM^*)=0$. Then $M$ is free.
\end{cor}

\begin{proof}
By \ref{lcrigid} and \ref{even} we have $\Tor_i^R(M,M^{*})=0$ for
all $i>0$. Thus Proposition \ref{3.2} shows that $M$ is free.
\end{proof}

\begin{eg}
Let $R,M$ as in Example \ref{mainex}. Then
$\LC_m^2(M\tensor_RM^*)=0$ but $M^*$ is not free. Note that $\dim
R=3$.
\end{eg}

\begin{cor}\label{lines}
Let $k$ be a field and $n$ an even integer $\geq 4$. Let $X \subset
\PP_k^{n}$ be a nonsingular hypersurface. Let $E$ be a vector bundle
on $X$. If $H^1(X, (E\tensor E^*)(l))=0$ for all $l \in \ZZ$, then
$E$ splits.
\end{cor}

\begin{proof}
This is a straightforward consequence of standard connections
between a projective variety and its affine cone (see section 5 in
\cite{HW2}). Let $A$ be the homogenous coordinate ring of $X$, $m$
be the irrelevant ideal and $R=A_m$. Let $N= \oplus_{i\in \ZZ
}H^0(X,E(i))$ and $M=N_m$. Note that the cohomology condition on $E$
translates to: $\LC_m^2(M\tensor_RM^*)=0$. Now we can apply
\ref{coreven2} to conclude that $M$ is a free $R$-module, which
means $E$ splits.
\end{proof}

It is worth noting  two obvious consequences of Corollary
\ref{coreven1} which are generalizations of well-known results by
Auslander and Goldman (\cite{AG}, Theorem 4.4) and Auslander
(\cite{Au}, Theorem 1.3) for modules over regular local rings:

\begin{cor}\label{coreven3}
Let $R$ be an admissible hypersurface such that $R$ has isolated
singularity and $\dim R>2$ and is even. Let $M$ be a reflexive
$R$-module. If $\Hom_R(M,M) \cong R^n$ for some $n$ then $M$ is
free.
\end{cor}

\begin{cor}\label{coreven4}
Let $R$ be an admissible hypersurface such that $R$ has isolated
singularity and $\dim R>2$ and is even. Let $M$ be an $R$-module
satisfying $(\Se_3)$. If $\Hom_R(M,M) \cong M^n$ for some $n$ then
$M$ is free.
\end{cor}

We now show that Conjecture \ref{conj2} is true in the standard
graded case:

\begin{thm}\label{projhyper}
Let $k$ be an algebraically closed field and  $n$ an odd integer.
Let $X \subset \PP_k^{n+1}$ be a smooth projective hypersurface. Let
$U,V$ be subvarieties of $X$ such that $U\cap V = \emptyset$.Then
$\dim U + \dim V < \dim X$.
\end{thm}

\begin{proof}
We are going to use $l$-adic cohomology (for basic properties and
notations, we refer to \cite{Mi} or \cite{Sha}). Let $l$ be a prime
number such that $l \neq char(k)$. There is a class map:
$$  cl : \CH^r(X) \to H^{2r}(X,\QQ_l (r))$$
This map gives a graded rings homeomorphism $\CH^*(X) \to \oplus
H^{2r}(X,\QQ_l (r))$ (with the intersection product on the left hand
side and the cup product on the right hand side, see \cite{Mi}, VI,
10.7 and 10.8). Let $a = \dim U$ and $b = \dim V$, and we may assume
$a\geq b$. Suppose $a+b = \dim X=n$ (if $a+b>n$, we can always
choose some subvariety of smaller dimension inside $U$ or $V$ such
that equality occurs). Then $2a \geq n$, but $n$ is odd, so $2a>n$.
Let $h \in \CH^1(X)$ represent the hyperplane section. By Weak
Lefschetz Theorem (see, for example, \cite{Sha}, 7.7, page 112) and
the fact that $2(n-a)<n$, we have:
$$H^{2(n-a)}(X,\QQ_l (n-a)) \cong H^{2(n-a)}(\PP_k^{N},\QQ_l
(n-a))$$ The latter is generated by a power of the class of the
hyperplane section. Thus $cl(U) = cl(h)^{n-a}$ in
$H^{2(n-a)}(X,\QQ_l (n-a))$. We then have :
$$ cl(U.V)= cl(U)\cup cl(V) = cl(h)^{n-a}\cup cl(V) =
cl(h^{n-a}.V)$$

The last term is equal to the degree of $h^{n-a}.V \in \CH^n(X)$, so
it is nonzero. But the first term has to be $0$ by assumption. This
contradiction proves the Theorem.
\end{proof}

\begin{cor}\label{graded}
Let $k$ be a perfect field and $n$ an even integer. Let  $A =
k[x_0,...,x_{n}]/(F)$ be a homogenous hypersurface. Let $(R,m,k)$ be
the local ring at the origin of $A$. Suppose that $R$ has an
isolated singularity. Let $M,N$ be graded $R$-modules such that
$l(M\tensor_RN)<\infty$. Then $\dim M +\dim N \leq \dim R$.
\end{cor}

\begin{proof}
Without affecting matters we may assume $k$ is algebraically closed
(it would not affect the isolated singularity condition because $k$
is perfect, so we can compute the singular locus of $A$ by the
Jacobian ideal, which would not be changed by extending $k$ to its
algebraic closure). Since the minimal primes of a graded modules are
homogenous, we may replace $M$ and $N$ with $R/P$ and $R/Q$, where
$P,Q$ are homogenous primes in $R$. Now let $X = \Proj(A)$,
$U=\Proj(A/P)$, $V=\Proj(A/Q)$ and applying the previous Theorem.
\end{proof}

\section{Some other applications}

In this section we discuss some further applications of rigidity.
The theme is that some strong conditions on a module could be
detected by the vanishing of a single $\Ext$ or $\Tor$ module. We
first note a simple characterization of maximal Cohen-Macaulay
modules, due to the fact that over an admissible hypersurface, a
module of finite length is rigid (see 2.4 in \cite{HW1} or 4.10 in
\cite{Da1}):
\begin{prop}\label{MCM}
Let $R$ be a admissible hypersurface and $M$ an $R$-module. The
following are equivalent:\\
(1) There is a nonzero finite length module $N$ such that $\Tor_1^R(M,N)=0$. \\
(2) There is a nonzero finite length module $N$ such that $\Ext_R^1(M,N)=0$. \\
(3) $M$ is a maximal Cohen-Macaulay $R$-module.
\end{prop}

\begin{proof}
For equivalence of (1) and (2): we may assume $R$ is complete. Then
we have a well-known isomorphism (see, for example, page 7 of
\cite{EG}) $\Tor_1^R(M,N)^{\vee} \cong \Ext_R^1(M,N^{\vee})$, so (1)
and (2) are equivalent.

Assume (1). Then since $N$ is rigid, we have $\Tor_i^R(M,N)=0$ for
all $i\geq 1$. The depth formula shows:
$$ \depth(M) + \depth(N) = \depth(R)+\depth(M\tensor_RN) $$
hence $\depth(M) = \depth(R)$ and (3) follows.

Finally, assume (3). Let $n=\dim R$, we can choose a regular
sequence $x_1,...,x_n$ on both $M$ and $R$ (by induction on $n$).
Now, just take $N = R/(x_1,...,x_n)$. Then $N$ has finite length and
$\Tor_1^R(M,N)=0$.
\end{proof}

It would be very interesting to know whether the previous
Proposition is true for all local complete intersections. It is not
hard to see that this question has intimate connection to the
rigidity of modules of finite length over such rings:

\begin{lem}
Let $\mathcal{C}$ be the category of all local complete
intersections. Consider the following properties:\\
(1) For any $R\in \mathcal{C}$, any $R$-module of finite length is
rigid.\\
(2) For any $R\in \mathcal{C}$ and any $R$-module $M$, if there is a
nonzero finite length $R$-module $N$ such that $\Tor_1^R(M,N)=0$
then $M$ is
maximal Cohen-Macaulay. \\
(3) For any $R\in \mathcal{C}$ such that $\dim R>0$ and any
$R$-module $M$, if there is a nonzero finite length $R$-module $N$
such that
$\Tor_1^R(M,N)=0$ then $\depth(M)>0$. \\
(4) For any $R\in \mathcal{C}$, any $R$-module of finite length and
finite projective dimension is rigid. \\
We have $(1) \Rightarrow (2) \Leftrightarrow (3) \Rightarrow (4)$.

\end{lem}

\begin{proof}
 $(1) \Rightarrow (2)$ is the proof of \ref{MCM}.

$(2) \Leftrightarrow (3)$. (2) clearly implies (3). Suppose (3)
holds, and let $M$ as in (2). We may assume $\dim R>1$. Then since
$\depth (M), \depth (R)>0$ we can find $x \in \Ann_R N$ such that
$x$ is regular on both $R$ and $M$. Then $\Tor_1^{R/(x)}(M/xM, N)
\cong \Tor_1^R(M,N)=0$. So we can replace $R,M$ by $R/(x), M/xM$ and
repeat the process until $\dim R=1$, at which point $M$ is MCM by
(3).

 $(2) \Rightarrow (4)$. Suppose $M,N$
are nonzero $R$-modules such that $l(N)<\infty$, $\pd_RN<\infty$ and
$\Tor_i^R(M,N)=0$ for some $i>0$. Then let $M' = \syz_{i-1}^RM$ and
we get $\Tor_1^R(M',N)=0$. By (2) $M'$ is MCM. Let $q$ be the
largest integer such that $\Tor_q^R(M',N) \neq 0$. Then by Lemma 2.2
in \cite{HW2} we have (since $\depth(\Tor_q^R(M',N))=0$):

$$ \depth (M') + \depth (N) = \depth (R) -q  $$
Since $\depth (M') = \depth (R)$ and $\depth (N) = 0$ we must have
$q=0$.
\end{proof}

\begin{rmk}
In \cite{JS} a module $M$ is constructed over a $0$-dimensional
Gorenstein ring such that $M$ is not rigid. So  property (1) fails
for Gorenstein rings.
\end{rmk}

Finally, we want to use our approach to rigidity to prove a
connection between vanishing of $\Ext$ and projective dimension,
generalizing a result by Jothilingam. In \cite{Jot}, the following
was proved:

\begin{thm}(\cite{Jot},Theorem)\label{Jor}
Let $R$ be a regular local ring and $M$  an $R$-module. Then
$\Ext^n_R(M,M)=0$ if and only if $\pd_R(M)<n$.
\end{thm}

In \cite{Jo}, Jorgensen observed that Jothilingam's result is true
for any local ring $R$ if we assume that  $M$ is rigid. We will
modify Jothilingam's proof to show:

\begin{prop}\label{Jo1}
Let $R$ be an admissible hypersurface and $M$ be an $R$-module such
that $M=0$ in $\overline G(R)_{\mathbb{Q}}$. Then $\Ext^n_R(M,M)=0$
if and only if $\pd_R(M)<n$.
\end{prop}

\begin{rmk}
Over a hypersurface, the classes of rigid modules and of modules
which are $0$ in $\overline G(R)_{\mathbb{Q}}$ are incomparable (see
\cite{Da1}).
\end{rmk}

\begin{proof}
We first need to review some notation in \cite{Jot}. Let $N$ be an
$R$-module and let
$$ F = ...\to F_2 \to F_1 \to F_0 \to N \to 0$$
be a minimal free resolution of $N$. Then one can define for $i\geq
0$:
$$ D^i(N) = \coker(F_i^* \to F_{i+1}^*)$$

Note that $D^i(N)$ can be computed, up to free summand, from any
resolution of $N$. We will use induction on $\dim R$. The proof of
the Theorem in \cite{Jot} shows that to prove $\Ext^n_R(M,M)=0$
implies $\pd_R(M)<n$ one only needs that the pair $(D^n(M), M)$ is
rigid. Let $L=D^n(M)$. We will show that $\theta^R(L,M)=0$. If $\dim
R=0$ then $\theta^R(-,-)$ is always defined and equals $0$, since
all modules have finite length and $\overline G(R)_{\mathbb{Q}}=0$.
Take any prime $p\in Y(R)$, the punctured spectrum. By induction
hypothesis $(\syz_{n-1}^RM)_p$ is free, so $L_p$ is a free $R_p$-
module. Thus $\theta^R(L,-)$ is always defined and since $M=0$ in
$\overline G(R)_{\mathbb{Q}}$ we must have $\theta^R(L,M)=0$. So
$(L,M)$ is rigid by \ref{rg1} and we are done.

\end{proof}

\begin{eg}
Let $R=k[[x,y]]/(xy)$ and $M =R/(x)$. Then $\Ext^{2i+1}(M,M)=0$ for
all $i>0$ but $\pd_RM=\infty$. Note that $\overline G(R)\cong \ZZ
[M] \cong \ZZ$.
\end{eg}

\section*{Ackowledgements}
We would like to thank Ragnar-Olaf Buchweitz, David Jorgensen, Craig
Huneke, Mircea Musta\c t\v a, Paul Roberts, Vasudevan Srinivas and
Claire Voisin for patiently explaining many useful facts and ideas
to us. We also thank the anonymous referee whose reports have
tremendously improved the presentation and mathematical content of
the note.

\end{document}